\documentclass[11pt,epsf]{article}
\usepackage{graphicx}
\usepackage{amsthm}
\usepackage{amsfonts}
\usepackage{amssymb}
\title{Stabilization of the asymptotic expansions 
of the zeros of a partial theta function}
\author{Vladimir Petrov Kostov\\ 
Universit\'e de Nice, 
Laboratoire de Math\'ematiques, Parc Valrose,\\ 06108 Nice Cedex 2, France,  
e-mail: vladimir.kostov@unice.fr} 
\date{}
\newtheorem{tm}{Theorem}

\newtheorem{rem}[tm]{Remark}

\begin{document} 
\maketitle 
\begin{abstract}
The bivariate series $\theta (q,x):=\sum _{j=0}^{\infty}q^{j(j+1)/2}x^j$ defines  
a {\em partial theta function}. For fixed $q$ ($|q|<1$), $\theta (q,.)$ is an 
entire function. We prove a property of stabilization of the 
coefficients of the Laurent series in $q$ of the zeros of 
$\theta$. These series are of the form 
$-q^{-j}+(-1)^jq^{j(j-1)/2}(1+\sum _{k=1}^{\infty}g_{j,k}q^k)$. 
The coefficients of the stabilized series 
are expressed by the positive integers $r_k$ giving the number 
of partitions into parts of three different kinds. They satisfy  
the recurrence relation 
$r_k=\sum _{\nu =1}^{\infty}(-1)^{\nu -1}(2\nu +1)r_{k-\nu (\nu +1)/2}$. 
Set $(H_{m,j})~:~(\sum _{k=0}^{\infty}r_kq^k)
(1-q^{j+1}+q^{2j+3}-\cdots +(-1)^{m-1}q^{(m-1)j+m(m-1)/2})=
\sum _{k=0}^{\infty}\tilde{r}_{k;m,j}q^k$. Then 
for $k\leq (m+2j)(m+1)/2-1-j$ and $j\geq (2m-1+\sqrt{8m^2+1})/2$ 
one has $g_{j,k}=\tilde{r}_{k;m,j}$.

{\bf AMS classification:} 26A06\\ 

{\bf Keywords:} partial theta function; asymptotics
\end{abstract}

The partial theta function is the sum of the bivariate series 
$\theta (q,x):=\sum _{j=0}^{\infty}q^{j(j+1)/2}x^j$, where 
$q\in \mathbb{C}$ ($|q|<1$) 
is considered as a parameter and $x\in \mathbb{C}$ as a variable. For 
each $q$ fixed, $\theta$ is an entire function. The terminology 
``partial theta function'' stems from the fact that, compared to the Jacobi 
theta function $\Theta (q,x):=\sum _{-\infty}^{\infty}q^{j^2}x^j$, 
only partial summation is 
performed (only over the non-negative indices), and 
$\theta (q,x)=\Theta (\sqrt{q},\sqrt{q}x)$.

The interest in the study of $\theta$ is explained by its applications 
in different domains. For instance, its relationship to hyperbolic polynomials 
(i.e. real polynomials with all roots real) has been exhibited in the 
recent papers \cite{Ost}, \cite{KaLoVi}, 
\cite{KoSh} and \cite{Ko2} which continue the earlier articles of 
Hardy, Petrovitch and Hutchinson \cite{Ha}, \cite{Hu} and \cite{Pe}. Other, 
better known domains, where $\theta$ is used, are the theory 
of (mock) modular forms (see \cite{BrFoRh}), 
asymptotic analysis (see \cite{BeKi}), statistical physics 
and combinatorics (see \cite{So}), Ramanujan type $q$-series 
(see \cite{Wa}); see also~\cite{AnBe}.

For $|q|\leq 0.108$ all zeros of $\theta$ are distinct, see~\cite{Ko4}. 
They can be indexed 
by the order of the pole at $0$ they have as functions of $q$. More precisely, 
the $j$-th zero of the partial theta function 
$\theta (q,x):=\sum _{j=0}^{\infty}q^{j(j+1)/2}x^j$ can be expanded in a 
Laurent series of the form 
$-\xi _j=-q^{-j}+(-1)^jq^{j(j-1)/2}(1+\sum _{k=1}^{\infty}g_{j,k}q^k)$, see~\cite{Ko6} 
(we set $g_{j,0}=1$). 
Also in \cite{Ko6} is shown that there exists a series 
$(H)~:~\sum _{k=0}^{\infty}r_kq^k$, where 
$r_k=\sum _{\nu =1}^{\infty}(-1)^{\nu -1}(2\nu +1)r_{k-\nu (\nu +1)/2}$, $r_0=1$, $r_k=0$ 
for $k<0$, with the property that 
for $k=1,\ldots ,j$ and $j\geq 2$ one has $g_{j,k}=r_k$. This property can 
be termed as stabilization of the Laurent series of the zeros $\xi _j$ 
as $j$ increases.  

In the present paper we improve this last property. 
We define the power series $(H_{m,j})$, $m,j\in \mathbb{N}$, by the formula:  
$$(H_{m,j})~:~(\sum _{k=0}^{\infty}r_kq^k)
(1-q^{j+1}+q^{2j+3}-\cdots +(-1)^{m-1}q^{(m-1)j+m(m-1)/2})=
\sum _{k=0}^{\infty}\tilde{r}_{k;m,j}q^k~.$$ 
Thus $(H_{1,j})=(H)$. We prove the following theorem:

\begin{tm}\label{TTTTT}
For $k\leq (m+2j)(m+1)/2-1-j$ and $j\geq (2m-1+\sqrt{8m^2+1})/2$ 
one has $g_{j,k}=\tilde{r}_{k;m,j}$.
\end{tm}

The sequence $\{ r_k\}$ is well-known, see \cite{Sl}. It gives the number of 
partitions into parts of three different kinds. 
We list its first $39$ elements:
$$1~, 3~, 9~, 22~, 51~, 108~, 221~, 429~, 810~, 1479~, 
2640~, 4599~, 7868~, 13209~, 21843~, 35581~, 57222~,$$ 
$$90882~, 142769~, 221910~, 
341649~, 521196~, 788460~, 1183221~, 1762462~, 2606604~,$$ 
$$3829437~, 5590110~, 
8111346~, 11701998~, 16790136~, 23964594~, 34034391~,$$ 
$$48104069~, 67679109~, 94800537~, 
132230021~, 183686994~, 254170332~.$$

{\em Proof of Theorem~\ref{TTTTT}:}\\ 

Set $-\xi _j=-q^{-j}+(-1)^jq^{j(j-1)/2}(\sum _{k=0}^{\infty}g_{j,k}q^k)$. For 
$\nu >1$ one has 
$$(-\xi _j)^{\nu}=(-q^{-j})^{\nu}+
\nu (-q^{-j})^{\nu -1}(-1)^jq^{j(j-1)/2}(\sum _{k=0}^{\infty}g_{j,k}q^k)+N~,$$ 
where $N$ contains all non-linear terms in $g_{j,k}$. Consider the series 
$\theta (q,-\xi _j)=\sum _{\nu =0}^{\infty}\Psi _{\nu}$, 
$\Psi _{\nu}=q^{\nu (\nu +1)/2}(-\xi _j)^{\nu}$. The lowest degree of $q$ in the 
expansion of $\Psi _{\nu}$ in Laurent series equals 
$\lambda _{\nu}:=\nu (\nu +1)/2-j\nu$. If $j\geq 2$, then the 
minimal degree of $q$ encountered in a coefficient of a non-linear term in 
the expansion of $\Psi _{\nu}$ is 
\begin{equation}\label{SSSSS}
j(j-1)-j(\nu -2)+\nu (\nu +1)/2=j(j+1-\nu )+\nu (\nu +1)/2~.
\end{equation} 
The coefficients $g_{j,k}$ can be computed from the condition 
$\theta (q,-\xi _j)=0$, by considering the coefficients of the powers of 
$q$ starting from the lowest one which is 
$\lambda _{j-1}=\lambda _j=\lambda _{\nu}=-j(j-1)/2$. 
When $j$ is sufficiently large compared to $k$, 
then non-linear terms do not intervene 
in these computations. In the proof of the theorem we explicit the conditions 
under which this does not take place.
 
The matrix that follows is denoted by $M_1$. Its columns contain the 
coefficients of the Laurent series in $q$ of $\Psi _0$, $\Psi _{j-3}$, $\ldots$, 
$\Psi _{j+2}$, $\Psi _{2j-1}$ and $\Psi _{2j}$. The first column indicates 
the power of $q$. For negative powers only the rows containing non-zero 
coefficients are represented. For brevity we set 
$a=g_{j,0}$, $b=g_{j,1}$, $c=g_{j,2}$, $\ldots$. The matrix 
corresponds to an index $j$ greater than $4$. (For $j=4$ the corresponding 
matrix is given in~\cite{Ko6}.) In the rows of higher powers of $q$ (which 
are not represented in the matrix) non-linear terms in $a$, $b$, $\ldots$ are 
present as well.   
$$\begin{array}{l|cccc|ccccc}
j+5&&(j-3)c&-(j-2)h&(j-1)u&-ju&(j+1)h&-(j+2)c&& \\ 
j+4&&(j-3)b&-(j-2)d&(j-1)h&-jh&(j+1)d&-(j+2)b&& \\ 
j+3&&(j-3)a&-(j-2)c&(j-1)d&-jd&(j+1)c&-(j+2)a&& \\ 
j+2&&&-(j-2)b&(j-1)c&-jc&(j+1)b&&& \\ 
j+1&&&-(j-2)a&(j-1)b&-jb&(j+1)a&&& \\ 
j&&&&(j-1)a&-ja&&&&1 \\ \hline  
1~{\rm to}~j-1&&&&&&&&& \\ \hline  
0&1&&&&&&&-1& \\ \hline 
\lambda _{j+2}&&(-1)^{j-3}&&&&&(-1)^{j+2}&& \\ 
\lambda _{j+2}-1&&&&&&&&& \\ 
\lambda _{j+1}&&&(-1)^{j-2}&&&(-1)^{j+1}&&& \\ 
\lambda _j&&&&(-1)^{j-1}&(-1)^j&&&& \\  
\hline  
&\Psi _0&\Psi _{j-3}&\Psi _{j-2}&\Psi _{j-1}&\Psi _j&\Psi _{j+1}&
\Psi _{j+2}&\Psi _{2j-1}&
\Psi _{2j} 
\end{array}$$
 
If one extends the matrix $M_1$ to the right, by adding the columns of 
$\Psi _{2j+1}$, $\Psi _{2j+2}$, $\ldots$, then the terms 
$(-1)^p$ appear in the columns of $\Psi _p$ and in the rows corresponding to 
$q^{p(p+1)/2-pj}$. Consider the restriction of the matrix $M_1$ 
to its rows not containing non-linear terms. Hence its $\mu$-th column 
(considered only in the rows corresponding to $q^s$ for $s\geq j$) 
up to a sign is of the form 
$(\ldots ,0,\mu a, \mu b, \mu c,\ldots )$. If there are no terms $(-1)^p$ 
for $p>2j$, the rows of the matrix give rise to the linear equations 
\begin{equation}\label{abc}
-a+1=0~,~-b+3a=0~,~-c+3b=0~,~-d+3c-5a=0~,~-h+3d-5b=0~~\ldots
\end{equation}
The solution to this system is the series $(H)$. 

Enlarge the matrix $M_1$ to the right by adding the column of $\Psi _{2j+1}$. 
This adds the term $-1$ in the column of $\Psi _{2j+1}$ and the row of 
$q^{2j+1}$ and no other non-zero terms.  
We denote by $M_2$ the new matrix thus obtained. For $\nu \leq j$ we set 
$g^0_{j,\nu}=g_{j,\nu}$ and for $\nu \geq j+1$ we set 
$g_{j,\nu}=g^0_{j,\nu}+g^*_{j,\nu}$, where $g^0_{j,\nu}$ are solutions to system 
(\ref{abc}), i.e. 
\begin{equation}\label{g0}
\begin{array}{rrr}-g^0_{j,0}+1=0~~&-g^0_{j,1}+3g^0_{j,0}=0~~&
-g^0_{j,2}+3g^0_{j,1}=0\\ \\ 
-g^0_{j,3}+3g^0_{j,2}-5g^0_{j,0}=0~~&-g^0_{j,4}+
3g^0_{j,3}-5g^0_{j,1}=0~~&\ldots ~~~~~\end{array}\end{equation}
and $g^*_{j,\nu}$ are solutions to the system
\begin{equation}\label{g*}
\begin{array}{rrr}-g^*_{j,j+1}-1=0~~&-g^*_{j,j+2}+3g^*_{j,j+1}=0~~&
-g^*_{j,j+3}+3g^*_{j,j+2}=0\\ \\ 
-g^*_{j,j+4}+3g^*_{j,j+3}-5g^*_{j,j+1}=0~~&-g^*_{j,j+5}+
3g^*_{j,j+4}-5g^*_{j,j+2}=0~~&\ldots ~~~~~\end{array}\end{equation}
The solution to system (\ref{g*}) is the series $-(H)$ in which the second 
indices of the unknown variables $g^*_{j,k}$ are shifted by $j+1$ compared to 
system (\ref{g0}). Hence 
the solution to the linear system resulting from the matrix 
$M_2$ is the series $(H_{2,j})=(H_{1,j})(1-q^{j+1})$. 

By adding one by one to the matrix the columns of 
$\Psi _{2j+2}$, $\Psi _{2j+3}$, $\ldots$ 
one obtains the matrices $M_3$, $M_4$, $\ldots$ which define 
linear systems whose solutions are the coefficients of the corresponding series 
$(H_{3,j})$, $(H_{4,j})$, $\ldots$. These solutions are the coefficients $g_{j,k}$ 
provided that the following two conditions are fulfilled:
\vspace{2mm}

(i) $k\leq (m+2j)(m+2j+1)/2-(m+2j)j-1-j=(m+2j)(m+1)/2-1-j$. Indeed,  
the last column of the matrix $M_m$ 
is the one of $\Psi _{m+2j-1}$. The index $k$ can take only these values for 
which the term $(-1)^p$ has not appeared in the row of $q^{p(p+1)/2-pj}$ for 
$p=m+2j$. These are the values for which the absence of the columns of 
$\Psi _{\nu}$ for $\nu \geq m+2j$ in $M_m$ does not affect the computation 
of the coefficients $g_{j,k}$). One has to subtract $j$ 
because the coefficient $a=g_{j,0}$ appears first in the row corresponding to 
$q^j$.    
\vspace{2mm}

(ii) The minimal power of $q$ multiplying a non-linear 
term is $\geq (m+2j)(m+2j+1)/2-(m+2j)j=(m+2j)(m+1)/2$ 
(hence the absence of the non-linear terms 
does not affect this computation either). This minimal power equals 

$$\min _{\nu}(j(j+1-\nu )+\nu (\nu +1)/2)=j(j+3)/2~~{\rm see~(\ref{SSSSS}).}$$ 
Hence one must have 

$$j(j+3)/2\geq (m+2j)(m+1)/2~,~~{\rm i.~e.}~~
j\geq (2m-1+\sqrt{8m^2+1})/2~.~~~~~\Box$$

\begin{rem}
{\rm In systems (\ref{abc}) and (\ref{g0}) the coefficients of the 
unknown variables are defined by the Laurent expansions of the monomials 
$\Psi _{\nu}$. Hence we implicitly suppose that no monomial $\Psi _{\nu}$ 
with $\nu <0$ 
(i.e. a nonexisting one) is involved. The coefficients 
$-1$, $3$, $-5$, $\ldots$  
of the variables $a$, $b$, $c$, $\ldots$ in these systems are obtained when 
considering the expansions of pairs of monomials 
$(\Psi _{\nu}, \Psi _{2j-1-\nu})$. 
The column of $\Psi _{\nu}$ (resp. $\Psi _{2j-1-\nu}$) 
gives the terms (see this column in the matrix) 
$\pm \nu a$, $\pm \nu b$, $\pm \nu c$, $\ldots$ 
(resp. $\pm (2j-1-\nu )a$, $\pm (2j-1-\nu )b$, 
$\pm (2j-1-\nu )c$, 
$\ldots$). Should there be involved a monomial $\Psi _{\nu}$ with $\nu <0$ in 
the computation of $a$, $b$, $c$, $\ldots$, then these nonexisting monomials 
should not be taken into account and the corresponding coefficients of systems 
(\ref{abc}) and (\ref{g0}) should be changed. 

The minimal positive powers of $q$ 
encountered in the pairs $(\Psi _1, \Psi _{2j-2})$, $(\Psi _0=1,\Psi _{2j-1})$ 
and $(\Psi _{-1},\Psi _{2j})$ (of which only the monomial $\Psi _{2j}$ exists) 
equal respectively $j(j-1)/2+1$, $j(j+1)/2$ and $j(j+3)/2$. The third 
of these numbers is 
precisely equal to the minimal power of $q$ multiplying a non-linear term, 
see (ii). For any monomial $\Psi _{\mu}$ with $\mu \geq 2j$ the minimal positive 
power of $q$ encountered in its Laurent expansion is $\geq j(j+3)/2$. 
Therefore pairs $(\Psi _{\nu}, \Psi _{2j-1-\nu})$ with negative values of $\nu$ 
are not involved in our reasoning.}
\end{rem}

\end{document}